\documentclass[final,leqno]{siamltex}
\usepackage{amsmath}
\usepackage{amssymb}
\usepackage{array}
\usepackage{algorithm}
\usepackage{algorithmic}
\usepackage{epsfig}

\usepackage{cite}

\usepackage{graphs}

\newcommand{\Real}{\mathbb R}
\newcommand{\Complex}{\mathbb C}
\newcommand{\set}[1]{\left\{#1\right\}}
\newcommand{\norm}[1]{\left\Vert#1\right\Vert}

\newcommand{\treeOne}{%
\begin{graph}(1.5,3)(-1,0)
\roundnode{rn1}(0,0)
\end{graph}}

\newcommand{\treeTwo}{%
\begin{graph}(1.5,3)(-1,0)
\roundnode{rn1}(0,0) \roundnode{rn2}(0,1) \edge{rn1}{rn2}
\end{graph}}

\newcommand{\treeThree}{%
\begin{graph}(1.5,3)(-1,0)
\roundnode{rn1}(0,0) \roundnode{rn2}(-.5,1) \roundnode{rn3}(.5,1)
\edge{rn1}{rn2} \edge{rn1}{rn3}
\end{graph}}

\newcommand{\treeFour}{%
\begin{graph}(1.5,3)(-1,0)
\roundnode{rn1}(0,0) \roundnode{rn2}(0,1) \roundnode{rn3}(0,2)
\edge{rn1}{rn2} \edge{rn2}{rn3}
\end{graph}}

\newcommand{\treeFive}{%
\begin{graph}(1.5,3)(-1,0)
\roundnode{rn1}(0,0) \roundnode{rn2}(-.75,1)
\roundnode{rn3}(.75,1) \roundnode{rn4}(0,1) \edge{rn1}{rn2}
\edge{rn1}{rn3} \edge{rn1}{rn4}
\end{graph}}

\newcommand{\treeSix}{%
\begin{graph}(1.5,3)(-1,0)
\roundnode{rn1}(0,0) \roundnode{rn2}(-.5,1) \roundnode{rn3}(.5,1)
\roundnode{rn4}(.5,2) \edge{rn1}{rn2} \edge{rn1}{rn3}
\edge{rn3}{rn4}
\end{graph}}

\newcommand{\treeSeven}{%
\begin{graph}(1.5,3)(-1,0)
\roundnode{rn1}(0,0) \roundnode{rn2}(0,1) \roundnode{rn3}(-.5,2)
\roundnode{rn4}(.5,2) \edge{rn1}{rn2} \edge{rn2}{rn3}
\edge{rn2}{rn4}
\end{graph}}

\newcommand{\treeEight}{%
\begin{graph}(1.5,3)(-1,0)
\roundnode{rn1}(0,0) \roundnode{rn2}(0,1) \roundnode{rn3}(0,2)
\roundnode{rn4}(0,3) \edge{rn1}{rn2} \edge{rn2}{rn3}
\edge{rn3}{rn4}
\end{graph}}

\title{A Family of ESDIRK Integration Methods
}

\author{John Bagterp J{\O}rgensen~\thanks{
Department of Applied Mathematics and Computer Science, Technical University of Denmark, DK-2800 Kgs.
Lyngby, Denmark ({\tt jbjo@dtu.dk})}\and Morten Rode
Kristensen\and {Per~Grove~Thomsen}}

\begin{document}

\maketitle

\begin{abstract}
In this paper we derive and analyze the properties of explicit
singly diagonal implicit Runge-Kutta (ESDIRK) integration methods.
We discuss the principles for construction of Runge-Kutta methods
with embedded methods of different order for error estimation and
continuous extensions for discrete event location. These
principles are used to derive a family of ESDIRK integration
methods with error estimators and continuous-extensions. The
orders of the advancing method (and error estimator) are 1(2),
2(3) and 3(4), respectively. These methods are suitable for
obtaining low to medium accuracy solutions of systems of ordinary
differential equations as well as index-1 differential algebraic
equations. The continuous extensions facilitates solution of
hybrid systems with discrete-events. Other ESDIRK methods due to
Kv{\ae}rn{\o} are equipped with continuous-extensions as well to
make them applicable to hybrid systems with discrete events.
\end{abstract}



\pagestyle{myheadings} \thispagestyle{plain} \markboth{J. B.
J{\O}RGENSEN, M. R. KRISTENSEN, P. G. THOMSEN}{A FAMILY OF ESDIRK
INTEGRATION METHODS}

\section{Introduction}
\label{sec:Introduction} In this paper, we derive and analyze a
family of explicit singly diagonally implicit Runge-Kutta (ESDIRK)
integration methods which can be applied for solution of stiff
systems of ordinary differential equations
\begin{equation}
\label{eq:ODEsystem}
    \dot{x}(t) = f(t,x(t)) \qquad x(t_0) = x_0
\end{equation}
as well as index-1 semi-explicit systems of differential algebraic
equations
\begin{subequations}
\label{eq:DAEsystem}
\begin{alignat}{3}
    \dot{x}(t) &= f(t,x(t),y(t)) \qquad && x(t_0) = x_0 \\
    0 &= g(t,x(t),y(t))
\end{alignat}
\end{subequations}
in which $t\geq t_0 \subset \Real$, $x \in \Real^n$, and $y \in
\Real^m$. For notational simplicity, we discuss the ESDIRK method
for (\ref{eq:ODEsystem}), but constructed with properties such
that it is equally applicable to (\ref{eq:DAEsystem}). A common
notation for (\ref{eq:ODEsystem}) and (\ref{eq:DAEsystem}) is
\begin{equation}
    M \dot{x}(t) = f(t,x(t)) \qquad x(t_0) = x_0
\end{equation}
in which the matrix $M$ may be singular. In addition it may depend
of $t$ and $x(t)$, i.e. $M = M(t,x(t))$.

Numerical methods for solution of these systems are not only of
importance in simulation, but are finding an increasing number of
applications in numerical tasks related to nonlinear predictive
control~\cite{Allgower:Zheng:2000,Betts:2001}. These tasks include
experimental design, parameter and state estimation, and numerical
solution of optimal control problems. While linear multistep
methods such as BDF based implementations, e.g.
DASPK~\cite{Gill:Jay:Leonard:Petzold:Sharma:2000,Serban:Petzold:2000,Maly:Petzold:1996},
DAEPACK~\cite{Park:Barton:1996,Feehery:Tolsma:Barton:1997,Tolsma:Barton:2000,Tolsma:Barton:2002,Barton:Lee:2002},
and
DAESOL~\cite{Bauer:Finocchi:Duschl:Gail:Schloder:1997,Bauer:Bock:Schloder:1999,Bauer:2000}
have been applied successfully to such problems, it has been
observed that typical industrial problems related to nonlinear
predictive control applications have frequent discontinuities.
Therefore, one-step methods, e.g.
SLIMEX~\cite{Schlegel:Marquardt:Ehrig:Nowak:2004} and ESDIRK
\cite{Kristensen:Jorgensen:Thomsen:Jorgensen:2004}, are more
efficient for the solution of such problems than linear multi-step
BDF methods.

Singly diagonally implicit Runge-Kutta methods (SDIRK) are
incepted by Butcher \cite{Butcher:1964} and have been applied for
solving systems of stiff ordinary differential equations since
their general introduction in the 1970s \cite{Norsett:1974,
Alexander:1977}. They have been provided in implementations such
as DIRKA and DIRKS~\cite{Cameron:1983, Cameron:Gani:1988}, SIMPLE
\cite{Norsett:Thomsen:1984,
Norsett:Thomsen:1986,Norsett:Thomsen:1986b}, and SDIRK4
\cite{Hairer:Wanner:1996}. SDIRK methods with an explicit first
stage equal to the last stage in the previous step are called
ESDIRK methods. They are of a much more recent origin and was
first considered as a general integration method for systems of
stiff systems in the years around year
2000~\cite{Cameron:1999,Cameron:Palmroth:Piche:2002,Butcher:Chen:2000,Williams:Cameron:Burrage:2000,Williams:Burrage:Cameron:Kerr:2002,Alexander:2003,Kvaernoe:2004}.
They retain the excellent stability properties of implicit
Runge-Kutta methods and compared to SDIRK methods improve the
computational efficiency. ESDIRK methods have been applied in
implicit-explicit Runge-Kutta methods for solution of
convection-diffusion-reaction
problems~\cite{Ascher:Ruuth:Spiteri:1997,Pareshci:Russo:2000,Kennedy:Carpenter:2003,Kristensen:Gerritsen:Thomsen:2006},
in dynamic optimization and optimal control applications for
efficient sensitivity
computation~\cite{Kristensen:Jorgensen:Thomsen:Jorgensen:2004,Kristensen:Jorgensen:Thomsen:Jorgensen:2004b,Kristensen:Jorgensen:IFAC:2005},
and for very computationally efficient implementations of extended
Kalman
filters~\cite{Jorgensen:Kristensen:Thomsen:Madsen:2006,Jorgensen:Kristensen:Thomsen:Madsen:2006b}.

In this paper, we derive and present a family of ESDIRK methods
suitable for numerical integration of stiff systems of
differential equations (\ref{eq:ODEsystem}) as well as index-1
differential algebraic systems (\ref{eq:DAEsystem}). The methods
are characterized in terms of A- and L-stability as well as order
of the basic integrator and the embedded method for error
estimation. We equip the methods with continuous-extensions such
that they can be applied to discrete-event systems. The paper is
organized as follows. In
Section~\ref{sec:RungeKuttaIntegrationMethods}, we present and
discuss Runge-Kutta methods and the general principles for their
implementation and construction. Section~\ref{sec:ESDIRKMethods}
applies these principles for construction of ESDIRK methods while
we discuss other ESDIRK methods in Sections
\ref{sec:ESDIRKKvaernoe} and \ref{sec:OtherESDIRKMethods}. The
other ESDIRK methods are equipped with continuous extensions.
Concluding remarks and a summarily comparison of the ESDIRK
methods are given in Section~\ref{sec:conclusion}. A companion
paper discusses implementation aspects and computational
properties of the ESDIRK
algorithms~\cite{Jorgensen:Kristensen:Thomsen:2006b}.

% //////////////////////////////////////////////////////////////////////////////

\section{Runge-Kutta Integration Methods}
\label{sec:RungeKuttaIntegrationMethods} The numerical solution of
systems of differential equations (\ref{eq:ODEsystem}) by an
s-stage Runge-Kutta method, may in each integration step be
denoted
\begin{subequations}
\label{eq:RungeKuttaMethodwithErrorEstimator}
\begin{alignat}{3}
    T_i &= t_n + c_i h \qquad && i=1,2, \ldots, s  \\
    X_i &= x_n + h \sum_{j=1}^{s} a_{ij} f(T_j,X_j) \qquad &&
    i=1,2, \ldots, s \label{eq:RungeKuttaMethodwithErrorEstimatorStageValues}\\
    x_{n+1} &= x_n + h \sum_{j=1}^{s} b_j f(T_j,X_j) \\
    \hat{x}_{n+1} &= x_n + h \sum_{j=1}^{s} \hat{b}_j f(T_j,X_j)
    \\
    e_{n+1} &= x_{n+1}-\hat{x}_{n+1}= h \sum_{j=1}^{s} d_j
    f(T_j,X_j) \quad && d_j = b_j - \hat{b}_j
\end{alignat}
\end{subequations}
$T_i$ and $X_i$ are the internal nodes and states computed by the
s-stage Runge-Kutta method. $x_{n+1}$ is the state computed at
$t_{n+1}=t_n+h$. $\hat{x}_{n+1}$ is the corresponding state
computed by the embedded Runge-Kutta method and $e_{n+1} = x_{n+1}
- \hat{x}_{n+1}$ is the estimated error of the numerical solution,
i.e. $\norm{e_{n+1}}$ is an estimate of the local error,
$\norm{x_{n+1} - x(t_{n+1})}$ given $x(t_n) = x_n$. The embedded
method, $\hat{x}_{n+1}$, uses the same internal stages as the
integration method, but the quadrature weights are selected such
that the embedded method is of different order, which then
provides an error estimate for the lowest order method. This order
relation of the integration method and the embedded method is
utilized by the error controller to adjust the step size, $h$,
adaptively~\cite{Gustafsson:1992,Gustafsson:1994,Soderlind:2002,Soderlind:2003}.

Alternatively, the s-stage Runge-Kutta method may be denoted and
implemented according to
\begin{subequations}
\label{eq:RungeKuttaMethodwithErrorEstimatorDerivative}
\begin{alignat}{3}
    T_i &= t_n + c_i h \qquad && i=1,2, \ldots, s \\
    X_i &= x_n + h \sum_{j=1}^{s} a_{ij} \dot{X}_j \qquad && i=1,2,
    \ldots, s \\
    \dot{X}_i &= f(T_i,X_i) = f(T_i, x_n + h \sum_{j=1}^{s} a_{ij}
    \dot{X}_j) \qquad && i=1,2, \ldots, s  \label{eq:RungeKuttaMethodwithErrorEstimatorDerivativeStageValue}\\
    x_{n+1} &= x_n + h \sum_{j=1}^{s} b_j \dot{X}_j \\
    \hat{x}_{n+1} &= x_n + h \sum_{j=1}^{s} \hat{b}_j \dot{X}_j \\
    e_{n+1} &= x_{n+1} - \hat{x}_{n+1} = h \sum_{j=1}^{s} d_j \dot{X}_j
    \quad &&d_j = b_j - \hat{b}_j
\end{alignat}
\end{subequations}
Sometimes the notation $K_i = \dot{X}_i$ is used for this
implementation. In (\ref{eq:RungeKuttaMethodwithErrorEstimator})
the stage values, $X_i$, are computed iteratively by solution of
(\ref{eq:RungeKuttaMethodwithErrorEstimatorStageValues}), while in
(\ref{eq:RungeKuttaMethodwithErrorEstimatorDerivative}) the time
derivatives of the stage values, $\dot{X}_i$, are computed
iteratively by solution of
(\ref{eq:RungeKuttaMethodwithErrorEstimatorDerivativeStageValue}).
Formally, (\ref{eq:RungeKuttaMethodwithErrorEstimator}) and
(\ref{eq:RungeKuttaMethodwithErrorEstimatorDerivative}) are
equivalent. However,
(\ref{eq:RungeKuttaMethodwithErrorEstimatorDerivative}) is
directly applicable to index-1 DAEs (\ref{eq:DAEsystem}) as well
as implicit DAE systems making this implementation preferred over
(\ref{eq:RungeKuttaMethodwithErrorEstimator}).

The s-stage Runge-Kutta method with an embedded error estimator,
(\ref{eq:RungeKuttaMethodwithErrorEstimator}) or
(\ref{eq:RungeKuttaMethodwithErrorEstimatorDerivative}), may be
denoted in terms of its Butcher tableau
\begin{equation*}
    \begin{array}{l|l}
        c & A \\ \hline
          & b' \\
          & \hat{b}' \\ \hline
          & d'
    \end{array}
    \qquad
    = \qquad
    \begin{array}{l|llll}
        c_1 & a_{11} & a_{12} & \ldots & a_{1s} \\
        c_2 & a_{21} & a_{22} & \ldots & a_{2s} \\
        \vdots & \vdots & \vdots & & \vdots \\
        c_s & a_{s1} & a_{s2} & \ldots & a_{ss} \\ \hline
        & b_1 & b_2 & \ldots & b_s \\
        & \hat{b}_1 & \hat{b}_2 & \ldots & \hat{b}_s \\ \hline
        & d_1 & d_2 & \ldots & d_s
    \end{array}
\end{equation*}
\begin{figure}[t]
\begin{center}
\epsfig{file=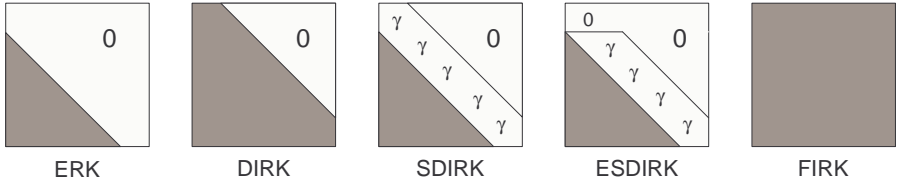,width=12cm} \caption{Structure of the
$A$-matrix for different classes of Runge-Kutta methods.}
\label{fig:RKAmatrix}
\end{center}
\end{figure}
Different classes of  Runge-Kutta methods may be characterized in
terms of the A-matrix in their Butcher tableau. This is
illustrated in Figure \ref{fig:RKAmatrix}. Explicit Runge-Kutta
(ERK) methods have a strictly lower triangular A-matrix implying
that (\ref{eq:RungeKuttaMethodwithErrorEstimatorDerivative}) may
be solved explicitly and without iterations. Therefore, ERK
methods have low computational cost but may suffer from stability
limitations when applied to stiff problems. ERK methods should
therefore be applied for non-stiff ODE problems, but not for stiff
ODE problems (\ref{eq:ODEsystem}) or DAE problems
(\ref{eq:DAEsystem}). All implicit Runge-Kutta methods are
characterized by an A-matrix that is not lower triangular. This
implies that some iterative method is needed for solution of
(\ref{eq:RungeKuttaMethodwithErrorEstimatorDerivative}). Fully
implicit Runge-Kutta (FIRK) methods are characterized by excellent
stability properties making them useful for solution of stiff
systems of ordinary differential equations (\ref{eq:ODEsystem}),
systems of index-1 semi-explicit differential algebraic equations
(\ref{eq:DAEsystem}), as well as systems of general differential
algebraic equations. However, in each integration step a system of
$n \times s$ coupled nonlinear equations must be solved. The price
of the excellent stability properties is high computational cost.
To achieve some of the stability properties of FIRK methods but at
lower computational cost, diagonally implicit Runge-Kutta (DIRK),
singly diagonally implicit Runge-Kutta (SDIRK), and explicit
singly diagonally implicit Runge-Kutta (ESDIRK) methods have been
constructed. For the DIRK methods, the internal stages decouples
in a such a way that the iterations may be conducted sequentially.
This implies that in DIRK-methods, $s$ systems of $n$ nonlinear
equations are solved instead of one system of $s\times n$
nonlinear equations as in the FIRK method. In the SDIRK methods,
the diagonal elements are identical such that the iteration matrix
may be reused for each stage. This saves a significant number of
LU-factorizations in the Newton iterations. In the ESDIRK method,
the first step is explicit ($c_1=0$ and $a_{11} = 0$), the
internal stages $2, \ldots, s$ are singly diagonally implicit, and
the last stage is equal to the next first stage ($c_s = 1$). This
implies that the first stage is free and that the iteration matrix
in stage $2, \ldots, s$ can be reused. Practical experience with
ESDIRK methods shows that they retain the stability properties of
FIRK methods but at significant lower computational costs.

In addition, ESDIRK methods are often constructed such that they
are {\em stiffly accurate}, i.e. $a_{si} = b_i$ for $i=1,2,
\ldots, s$ (note $b_s = a_{ss} = \gamma$). This implies that the
last stage is equal to the final solution, $x_{n+1} = X_s$, and
that no extra computations are needed for solution of the
algebraic variables in (\ref{eq:DAEsystem}). Furthermore, stiffly
accurate methods avoid the order reduction for stiff
systems~\cite{Proherto:Robinson:1974}. Stiffly accurate ESDIRK
methods with different number of stages and order can be
represented by the following Butcher tableaus
\begin{equation*}
\begin{array}{l|ll}
    0 & 0 & \\
    1 & b_{1} & \gamma \\ \hline
      & b_{1} & \gamma \\
      & \hat{b}_1 & \hat{b}_2 \\ \hline
      & d_1 & d_2
\end{array}
\qquad \quad
\begin{array}{l|lll}
    0 & 0 & & \\
    c_2 & a_{21} & \gamma & \\
    1   & b_1 & b_2 & \gamma \\ \hline
        & b_1 & b_2 & \gamma \\
        & \hat{b}_1 & \hat{b}_2 & \hat{b}_3 \\ \hline
        & d_1 & d_2 & d_3
\end{array}
\qquad \quad
\begin{array}{l|llll}
    0   & 0 & & & \\
    c_2 & a_{21} & \gamma & & \\
    c_3 & a_{31} & a_{32} & \gamma & \\
    1   & b_1 & b_2 & b_3 & \gamma \\ \hline
        & b_1 & b_2 & b_3 & \gamma \\
        & \hat{b}_1 & \hat{b}_2 & \hat{b}_3 & \hat{b}_4 \\ \hline
        & d_1 & d_2 & d_3 & d_4
\end{array}
\end{equation*}
As is evident from the above Butcher tableaus, stiffly accurate
ESDIRK methods with $s$ stages need only to compute $s-1$ stages
as the first stage is equal to the last stage of the previous
step.

\subsection{Order Conditions for Runge-Kutta Methods}
The order conditions for Runge-Kutta methods are developed
considering Taylor expansions of the analytical and numerical
solution of the autonomous ODE
\begin{equation}
\label{eq:AutonomousODEsystem}
    \dot{x}(t) = f(x(t)) \qquad x(t_0) = x_0
\end{equation}
The forced ODE (\ref{eq:ODEsystem}) can always be transformed into
an autonomous system of ODEs (\ref{eq:AutonomousODEsystem}) using
\begin{equation*}
    z(t) = \begin{bmatrix} x(t) \\ t \end{bmatrix} \quad
    \dot{z}(t) = \begin{bmatrix} \dot{x}(t) \\ \dot{t}
    \end{bmatrix} = \begin{bmatrix} f(t,x(t)) \\ 1 \end{bmatrix} =
    F(z(t)) \quad z(t_0) = \begin{bmatrix} x(t_0) \\ t_0
    \end{bmatrix} = \begin{bmatrix} x_0 \\ t_0 \end{bmatrix} = z_0
\end{equation*}
Runge-Kutta methods must satisfy the consistency condition
(\ref{eq:ConsistencyConditionRK}) to integrate the time component
correctly, i.e. to have exact equivalence between the forced
system (\ref{eq:ODEsystem}) and the autonomous system
(\ref{eq:AutonomousODEsystem}). The order conditions are usually
derived considering a scalar autonomous system. This is simpler
and loses no generality compared to the vector case. Let
$\mathcal{T}$ denote all rooted trees and let $\mathcal{T}(p)$
denote the set of rooted trees with order less or equal to $p$,
i.e. $\mathcal{T}(p) = \set{\tau \in \mathcal{T}: \, r(\tau) \leq
p}$. Some functions on rooted trees are listed in Table
\ref{table:FunctionsOnRootedTrees}. We call these trees for
Butcher trees since they were first used by Butcher to derive
order conditions for Runge-Kutta
methods~\cite{Butcher:1963,Butcher:Wanner:1996,Butcher:2003}. The
advantage of Butcher trees is that the order conditions can be
derived from these trees and it is significantly easier to write
up all Butcher trees to a given order than {\em ab initio}
derivation of the order conditions from Taylor expansions. The
number of nodes (dots) in the tree corresponds to the order,
$r(\tau)$, and the symmetry, $\sigma(\tau)$, is easily inspected
for a tree by labelling the nodes. The density, $\gamma(\tau)$, is
computed by multiplying the orders of each subtree rooted on a
vertex of $\tau$, e.g. $\gamma(\tau_7) = 4 \cdot  3 \cdot (1 \cdot
1) = 12$. The matrices, $\Lambda(\tau)$, can also be derived by
inspection of the Butcher trees. $\Lambda(\tau)$ is constructed as
follows: a vertex connecting to a node with no further subtrees
corresponds to multiplying by $C$, while a vertex connecting to a
node with further subtrees corresponds to multiplying by $A$. As
an example consider the tree $\tau_6$. The first vertex going to
the left ends on a terminal node and therefore corresponds to
multiplying by $C$. The first vertex going to the right does not
end on a terminal node and corresponds to multiplying by $A$,
while the next vertex on this branch of the tree ends on a
terminal node and corresponds to multiplying with $C$. Hence,
$\Lambda(\tau_6) = C A C$. $\Phi(\tau)$ and $\Psi(\tau)$ are
defined as $\Phi(\tau) = b' \Lambda(\tau) e$ and $\Psi(\tau) =
\Lambda(\tau) e$. The elementary weights $F(\tau)$ are also easily
derived from the rooted tree (see Table
\ref{table:FunctionsOnRootedTrees}).

\begin{table}[tb]
\caption{Some functions on rooted
trees~\cite{Butcher:Wanner:1996,Butcher:2003}. $\tau$ denotes the
rooted tree, $r(\tau)$ the order, $\sigma(\tau)$ the number of
symmetries, and $\gamma(\tau)$ the density computed by multiplying
the orders of each subtree rooted on a vertex of $\tau$.
$\Lambda(\tau)$ are elementary weights derived under the
consistency assumption $Ae = Ce$. $\Phi(\tau)$ and $\Psi(\tau)$
are defined by $\Phi(\tau) = b' \Lambda(\tau) e$ and $\Psi(\tau) =
\Lambda(\tau)e$. $F(\tau)$ denotes an elementary weight of $\tau$,
e.g. of its use: $F(\tau_3)(x(t_n)) = f''(f,f)(x(t_n)) =
f''(x(t_n)) f(x(t_n)) f(x(t_n))$.
}\label{table:FunctionsOnRootedTrees} \footnotesize
\begin{center}
\begin{tabular}{lcccccccc} \hline
 $\tau$ & $\tau_1$ & $\tau_2$ & $\tau_3$ & $\tau_4$  & $\tau_5$ &  $\tau_6$ & $\tau_7$ & $\tau_8$ \\
  & \treeOne & \treeTwo & \treeThree & \treeFour & \treeFive & \treeSix & \treeSeven & \treeEight \\ \hline
 $r(\tau)$ & 1 & 2 & 3 & 3 & 4 & 4 & 4 & 4 \\
 $\sigma(\tau)$ & 1 & 1 & 2 & 1 & 6 & 1 & 2 & 1 \\
 $\gamma(\tau)$ & 1 & 2 & 3 & 6 & 4 & 8 & 12 & 24 \\
 $\Lambda(\tau)$ & $I$ & $C$ & $C^2$ & $AC$ & $C^3$ & $CAC$ &
 $AC^2$ & $A^2C$ \\
 $\Phi(\tau)$ & $b' e$ & $b' C e$ & $b' C^2 e$ & $b' A C e$ & $b' C^3 e$ & $b' C A C e$ & $b' A C^2 e$ & $b' A^2 C e$  \\
 $\Psi(\tau)$ & $e$ & $C e$ & $C^2 e$ & $A C e$ & $C^3 e$ & $C A C e$ & $A C^2 e$ & $A^2 C e$  \\
 $F(\tau)$ & $f$ & $f'f$ & $f'' (f,f)$ & $f' f' f$ & $f'''(f,f,f)$ & $f''(f, f'
  f)$ & $f' f''(f,f)$ & $f' f' f' f$
 \\ \hline
\end{tabular}
\end{center}
\end{table}
Using Butcher trees and assuming $x(t_n) = x_n$,  the order
conditions for Runge-Kutta methods are derived by comparing the
Taylor expansion of the exact solution
\begin{equation}
\begin{split}
    x(t_{n+1}) = x(t_n+h) &= x(t_n) + \sum_{k=1}^{p} \frac{1}{k!} \frac{d^k
    x}{dt^k}(t_n) h^k + O(h^{p+1})
    \\&= x(t_n) + \sum_{\tau \in \mathcal{T}(p)}
    \frac{h^{r(\tau)}}{\sigma(\tau) \gamma(\tau)} F(\tau)(x(t_n))
    + O(h^{p+1})
\end{split}
\end{equation}
and the Taylor expansion of the numerical solution
\begin{equation}
    x_{n+1} = x_n + \sum_{\tau \in \mathcal{T}(p)}
    \frac{\Phi(\tau) h^{r(\tau)}}{\sigma(\tau)} F(\tau)(x_n) +
    O(h^{p+1})
\end{equation}
which is obtained by Taylor expansion of $f(X_i)$ in
(\ref{eq:RungeKuttaMethodwithErrorEstimator}) around $x_n$. The
method has order $p$ if the local error is $e_{n+1} = x_{n+1} -
x(t_{n+1}) = O(h^{p+1})$, i.e. if $\Phi(\tau) = 1/\gamma(\tau), \,
\forall \tau \in \mathcal{T}(p)$.

Let $C = \text{diag}\set{c_1, c_2, \ldots, c_s}$ be a diagonal
matrix with $\set{c_i}_{i=1}^{s}$ on the diagonal. Then the
consistency condition (the row-sum condition) can be expressed as
\begin{subequations}
\label{eq:OrderConditionsRK}
\begin{equation}
\label{eq:ConsistencyConditionRK}
    C e = A e
\end{equation}
and the order conditions $\Phi(\tau) = b' \Psi(\tau) =
1/\gamma(\tau)$ for $r(\tau)\leq p$ for order $p=\set{1,2,3,4}$
can be expressed as
\begin{alignat}{3}
& \text{Order } 1: \quad && b' e = 1 \\
& \text{Order } 2: \quad && b' C e = \frac{1}{2} \\
& \text{Order } 3: \quad && b' C^2 e = \frac{1}{3} \\
& && b' A C e = \frac{1}{6} \\
& \text{Order } 4: \quad && b' C^3 e = \frac{1}{4} \\
& && b' C A C e = \frac{1}{8} \\
& && b' A C^2 e = \frac{1}{12} \\
& && b' A^2 C e = \frac{1}{24}
\end{alignat}
\end{subequations}
Using the order conditions (\ref{eq:OrderConditionsRK}), the
special structure of the Butcher tableau of ESDIRK methods, and
the A- and L-stability conditions, we can derive ESDIRK methods of
various order. It turns out, that these conditions do not always
determine the methods uniquely. Therefore, we consider the
simplifying assumptions of Runge-Kutta methods as additional
design
criteria~\cite{Butcher:1963,Hairer:Norsett:Wanner:1993,Butcher:2003}
\begin{subequations}
\begin{alignat}{5}
& B(q): \qquad && b' C^{k-1} e = \frac{1}{k} \qquad && k=1,2,
\ldots, q
\\
& C(q): \qquad && A C^{k-1} e = \frac{1}{k} C^k e \qquad && k=1,2,
\ldots, q
\\
& D(q): \qquad && A' C^{k-1} b = \frac{1}{k} (I- C^k) b \qquad &&
k = 1,2, \ldots, q
\end{alignat}
\end{subequations}
The conditions $B(p)$ are part of the conditions for order $p$ and
therefore not included twice. We disregard the conditions $D(q)$.
This leaves the conditions $C(q)$. We note that $C(1)$ corresponds
to the consistency conditions (\ref{eq:ConsistencyConditionRK})
and $k=1$ is not included in $C(q)$. $C(q)$ implies that the
internal stages has order $q$, i.e. $E_i = X_i - x(t_n+c_ih) =
O(h^{q+1})$ \cite{Hairer:Wanner:1996}. For ESDIRK methods, $c_1 =
0$ and $c_s = 1$. The order conditions ensure that numerical
solutions at these points are of order $p$. Therefore, we do not
enforce $C(q)$ at these points and this leaves the following
additional design conditions
\begin{alignat}{5}
& \tilde{C}(q): \qquad && \sum_{j=1}^{s} a_{ij} c_j^{k-1} =
\frac{1}{k}c_i^k \qquad && i=2, 3,  \ldots, s-1; \, k=2, 3,
\ldots, q
\end{alignat}
implying stage order $q$ for the stages $i=2,3, \ldots, s-1$. The
methods considered in this paper has stage order 2, i.e. they
satisfy
\begin{equation}
\label{eq:ModifiedStageOrder2} \tilde{C}(2): \qquad \sum_{j=1}^{s}
a_{ij} c_j = \frac{1}{2} c_i^2 \qquad i = 2, 3,
    \ldots, s-1
\end{equation}
In particular, for ESDIRK methods, $\tilde{C}(2)$ implies that the
second stage, $i=2$, satisfies
\begin{equation}
    \gamma c_2 = \frac{1}{2} c_2^2 \quad \Leftrightarrow \quad
    c_2 = 2 \gamma
\end{equation}
Along with having A- and L-stability, the methods considered in
this paper satisfy the consistency and order conditions
(\ref{eq:OrderConditionsRK}) as well as stage order 2 conditions
(\ref{eq:ModifiedStageOrder2}).

\subsection{Continuous Extension}
\label{sec:ContinuousExtension} The ability to efficiently compute
a numerical approximation $\bar{x}(t_n+\theta h)$ to $x(t_n+\theta
h)$ for $\theta \in [0,1]$ is important for hybrid systems with
discrete events as well as in creating dense outputs for plotting
and visualization purposes. $\bar{x}(t_n+\theta h)$ is called a
continuous extension of the Runge-Kutta method. It is computed as
\begin{equation}
    \bar{x}(t_n+\theta h) = x_n + h \sum_{i=1}^{s} \bar{b}_i(\theta)
    \dot{X}_i
\end{equation}
in which
\begin{equation}
    \bar{b}(\theta) = \sum_{k=1}^{q} \bar{b}_k \theta^k \qquad
    \bar{b}(\theta) = \begin{bmatrix} \bar{b}_1(\theta) \\ \vdots
    \\ \bar{b}_s(\theta) \end{bmatrix} \quad \bar{b}_{k} =
    \begin{bmatrix} \bar{b}_{1k} \\ \vdots \\ \bar{b}_{sk}
    \end{bmatrix} \, k=1,2,\ldots, q
\end{equation}
The continuous extension is of order $q$ if $\bar{e}(t_n+\theta h)
= \bar{x}(t_n+\theta h) - x(t_n + \theta h) = O(h^{q+1})$ for
$x(t_n) = x_n$. To construct a continuous extension of order $q$,
we determine the coefficient matrix, $\bar{B} =
[\bar{b}_k]_{k=1,2, \ldots,q}$, such that it satisfies the
Runge-Kutta order conditions
\begin{equation}
    \forall \theta \in [0,1]: \quad \bar{b}(\theta)' \Psi(\tau) =
    \frac{\theta^{r(\tau)}}{\gamma(\tau)} \qquad \forall \tau \in
    \mathcal{T}(p)
\end{equation}
If the order, $q$, of the continuous extension is equal to the
order, $p$, of the advancing integration method we require
$\bar{x}(t_n+h) = x_{n+1}$ which corresponds to the condition
\begin{subequations}
\label{eq:RkContinuousExtensionLinearSideCondition}
\begin{equation}
    \bar{b}(\theta=1) = \sum_{k=1}^{p} \bar{b}_k =  b
\end{equation}
Similarly, if the order, $q$, of the continuous extension is equal
to the order, $\hat{p}$, of the embedded method we require
$\bar{x}(t_n+h) = \hat{x}_{n+1}$ which corresponds to the
condition
\begin{equation}
    \bar{b}(\theta=1) = \sum_{k=1}^{\hat{p}} \bar{b}_k = \hat{b}
\end{equation}
\end{subequations}
Consequently, the coefficients, $\bar{B} =
[\bar{b}_k]_{k=1,2,3,4}$, for the continuous extension of order
1-4 may be obtained as the solution of a linear system in the
following form
\begin{subequations}
\renewcommand{\arraystretch}{1.2}
\label{eq:RKContinuousExtensionLinearEquationsI}
\begin{gather}
    \overbrace{\left[
    \begin{array}{c}
        e' \\ \hline
        (C e)' \\ \hline
        (C^2 e)' \\
        (A C e)' \\ \hline
        (C^3 e)' \\
        (CAC e)' \\
        (A C^2 e)' \\
        (A^2 C e)'
    \end{array}
    \right]}^{\bar{\Psi} = }
    \overbrace{\left[
        \begin{array}{c|c|c|c}
            \bar{b}_1 & \bar{b}_2 & \bar{b}_3 & \bar{b}_4
        \end{array}
    \right]}^{\bar{B} =}
    =
    \overbrace{
    \left[
    \begin{array}{c|c|c|c}
        1 & 0 & 0 & 0 \\ \hline
        0 & \frac{1}{2} & 0 & 0 \\ \hline
        0 & 0 & \frac{1}{3} & 0 \\
        0 & 0 & \frac{1}{6} & 0 \\ \hline
        0 & 0 & 0 & \frac{1}{4} \\
        0 & 0 & 0 & \frac{1}{8} \\
        0 & 0 & 0 & \frac{1}{12} \\
        0 & 0 & 0 & \frac{1}{24}
    \end{array}
    \right]}^{\bar{\Gamma}=}
\\
    \left[
    \begin{array}{c|c|c|c} \bar{b}_1 & \bar{b}_2 & \bar{b}_3 & \bar{b}_4
    \end{array} \right] e = b \, \text{(or $\hat{b}$)}
\end{gather}
\end{subequations}
and the solution computed by solving
\begin{equation}
\label{eq:ContinuousExtensionOrderAndSideConditionsLinear}
    \begin{bmatrix}
        I_q \otimes \bar{\Psi} \\
        e_q' \otimes I_s
    \end{bmatrix}
    \text{vec}(\bar{B})
    = \begin{bmatrix} \text{vec}(\bar{\Gamma}) \\ b \,\text{(or $\hat{b}$)} \end{bmatrix}
\end{equation}
$\otimes$ denotes the Kronecker product, $\text{vec}$ denotes
vectorization of a matrix, $I_q$ is a q-dimensional unity matrix,
and $e_q' = \begin{bmatrix}1 & \ldots & 1\end{bmatrix} \in
\Real^q$. We have used the relation $\text{vec}(ABC) = (C' \otimes
A)\text{vec}(B)$ in the derivation of
(\ref{eq:ContinuousExtensionOrderAndSideConditionsLinear}). The
horizontal and vertical lines in
(\ref{eq:RKContinuousExtensionLinearEquationsI}) indicates which
parts to retain for various order of the continuous extension.

Other conditions may be considered as well for construction of the
continuous extension, i.e. $\bar{x}(t_n + c_i h) = X_i$ which
leads to a condition of the type
\begin{equation}
    \left[
        \begin{array}{c|c|c|c}
            \bar{b}_1 & \bar{b}_2 & \bar{b}_3 & \bar{b}_4
        \end{array}
    \right]
    \begin{bmatrix} c_i \\ c_i^2 \\ c_i^3 \\ c_i^4 \end{bmatrix}
    =
    \begin{bmatrix} a_{i1} \\ a_{i2} \\ a_{i3} \\ a_{i4}
    \end{bmatrix}
\end{equation}
and this may be incorporated in the linear system in a similar way
to the incorporation of
(\ref{eq:RkContinuousExtensionLinearSideCondition}). A derivative
condition, $\dot{\bar{x}}(t_n+c_ih) = \dot{X}_i$, requires
satisfaction of a linear constraint of the type
\begin{equation}
    \left[
    \begin{array}{c|c|c|c}
        \bar{b}_1 & \bar{b}_2 & \bar{b}_3 & \bar{b}_4
    \end{array}
    \right]
    \begin{bmatrix} 1 c_i^0 \\ 2 c_i^1 \\ 3 c_i^2 \\ 4 c_i^3
    \end{bmatrix}
    =
    e_i
\end{equation}
in which $e_i = \begin{bmatrix} 0 & \ldots & 0 & 1 & 0 & \ldots &
0\end{bmatrix}' \in \Real^s$ with the unit entry in the $ith$
coordinate.

\subsection{Stability Conditions}
Using the test equation $\dot{x}(t) = \lambda x(t)$ with initial
condition $x(0)=x_0$ and $\lambda \in \Complex$, all Runge-Kutta
methods (\ref{eq:RungeKuttaMethodwithErrorEstimatorDerivative})
for advancing the solution may be expressed as $x_{n+1} =
R(h\lambda) x_n$ in which the transfer function $R(z)$ is
\begin{equation}
    R(z) = 1 + z b' (I- z A)^{-1} e = \frac{\det\left( I - zA + z e b' \right)}{\det\left(I - zA \right)} = \frac{P(z)}{Q(z)}  \qquad z \in \Complex
\end{equation}
with $I$ being the unity matrix and $e = \begin{bmatrix} 1 & 1 &
\ldots & 1 \end{bmatrix}'$. An integration method is said to be
{\em $A$-stable} if its transfer function $R(z)$ for the test
equation is stable in the left half plane, i.e. if $|R(z)| < 1$
for $\text{Re}(z) < 0$. This implies that for Lyapunov stable test
equations, i.e. $\text{Re}(\lambda) < 0$, the numerical solution,
${x_n = R(h \lambda)^n x_0}$, obtained by the integration method
will converge to the mathematical solution, $x(t_n) = \left( e^{h
\lambda } \right)^n x_0$ with $t_n = nh$.

An integration method is said to be {\em $L$-stable} if its
transfer function, $R(z)$, for the test equation is A-stable and
in addition satisfies
\begin{equation}
    \lim_{z\rightarrow -\infty} |R(z)| = 0
\end{equation}
L-stability is an important property, when the integration method
is applied for solution of systems of differential algebraic
equations. Note that $|R(-\infty)| = |R(\infty)|$ for Runge-Kutta
methods. Consider a stiffly accurate ESDIRK method, i.e. a method
with the Butcher tableau
\begin{equation*}
\renewcommand{\arraystretch}{1.2}
    \begin{array}{l|l}
        c & A \\ \hline
         & b'
    \end{array}
    \quad = \quad
    \begin{array}{l|l|lllll}
        0 & 0 & 0 & 0 & \ldots & 0 & 0\\ \hline
        c_2 & a_{21} & \gamma & 0 & \ldots  & 0 & 0\\
        c_3 & a_{31} & a_{32} & \gamma & \ldots & 0 & 0  \\
        \vdots & \vdots & \vdots & \vdots & & \vdots & \vdots \\
        c_{s-1} & a_{s-1,1} & a_{s-1,2} & a_{s-1,3} & \ldots & \gamma & 0\\
        1 & b_1 & b_2 & b_3 & \ldots & b_{s-1} & \gamma \\ \hline
          & b_1 & b_2 & b_3 & \ldots & b_{s-1} & \gamma
    \end{array}
    \quad = \quad
    \begin{array}{l|ll}
        0 & 0 & 0 \\
        \tilde{c} & \tilde{a} & \tilde{A} \\ \hline
                  & b_1 & \tilde{b}'
    \end{array}
\end{equation*}
then~\cite{Kvaernoe:2004,Jay:1993}
\begin{equation}
    R(\infty) = - e_{s-1}' \tilde{A}^{-1} \tilde{a} \qquad
    e_{s-1}'
    = \begin{bmatrix} 0 & 0 & \ldots & 0 & 1 \end{bmatrix}
\end{equation}

For stiffly accurate s-stage ESDIRK methods, the numerator
polynomial $P(z) = \det\left(I - zA + z e b'\right)$ is at most of
degree $s-1$ and given by
\begin{equation}
    P(z) = (-1)^{s-1} \sum_{j=0}^{s-1}
    L_{s-1}^{(s-1-j)}\left(\frac{1}{\gamma}\right)(\gamma z)^j
\end{equation}
in which
\begin{equation}
    L_{s-1}(x) = \sum_{j=0}^{s-1} (-1)^j \begin{pmatrix} s-1 \\ j
    \end{pmatrix} \frac{x^j}{j!}
\end{equation}
are the Laguerre-polynomials and $L_s^{(k)}(x)$ denotes their
$k$th derivative. For s-stage ESDIRK methods, the denominator
polynomial is
\begin{equation}
Q(z) = \det\left(I - zA\right) = (1-\gamma z)^{s-1}
\end{equation}
Since $Q(z)$ is of degree $s-1$, the requirement of L-stability
corresponds to a zero coefficient for the term $z^{s-1}$ in the
numerator polynomial, i.e.
\begin{equation}
    \gamma \neq 0 : \quad
     (-1)^{s-1} L_{s-1}^{}\left(\frac{1}{\gamma}\right) \gamma^{s-1}
     =0 \quad \Leftrightarrow \quad
     L_{s-1}\left(\frac{1}{\gamma}\right) = 0
\end{equation}

The stability function of a stiffly accurate s-stage ESDIRK method
is identical to the stability function of an (s-1)-stage SDIRK
method~\cite{Norsett:1974,Hairer:Wanner:1996,Kvaernoe:2004}.
Hairer and Wanner \cite{Hairer:Wanner:1996} provide regions of A-
and L-stability of SDIRK methods. These regions and conditions are
translated into conditions for stiffly accurate ESDIRK methods and
listed in Table \ref{table:StabilityRegionsESDIRK}.

The last column in Table \ref{table:StabilityRegionsESDIRK}
indicates the location of the second quadrature point, $T_2 = t_n
+ h c_2$, provided stage order 2 is required for the second step,
i.e. $c_2 = 2 \gamma$ for order $p = s-1 \geq 2$. For one-step
methods it is reasonable for computational and implementation
reasons to require the quadrature points to be within the step,
i.e. $t_n \leq T_i \leq t_n + h$ which implies $0 \leq c_2 \leq 1$
or $0 \leq \gamma \leq \frac{1}{2}$. From Table
\ref{table:StabilityRegionsESDIRK} it is apparent that this
condition along with the requirements of A- and L-stability imply
that s-stage ESDIRK methods with order $p=s-1$ exist for
$s=\set{2,3,4}$ but not for $s=5$. In the cases $s=\set{2,3,4}$,
the requirements determine $\gamma$ uniquely.

\begin{table}[tb]
\caption{Stability of stiffly accurate s-stage ESDIRK methods of
order $p=s-1$~\cite[pp. 96-98]{Hairer:Wanner:1996}.}
\label{table:StabilityRegionsESDIRK} \footnotesize
\begin{center}
\begin{tabular}{llll}
    \hline
    s & A-stability & L-stability  &   \\
      & $p\geq s-1$ &              $p=s-1$ &
    \\ \hline
    2 & $\frac{1}{2} \leq \gamma < \infty$ & $\gamma = 1$ & $c_2 = 1$ \\
    3 & $\frac{1}{4} \leq \gamma < \infty$  & $\gamma = \frac{2 \pm \sqrt{2}}{2} = \begin{cases} 1.70710678 \\ 0.29289322 \end{cases}$ & $c_2 =
    2\gamma = \begin{cases} 3.41421356 \\ 0.58578644 \end{cases}$ \\
    4 & $\frac{1}{3} \leq \gamma \leq 1.06857902$  & $\gamma = 0.43586652 $ & $c_2 = 2 \gamma = 0.87173304 $ \\
    5 & $0.39433757 \leq \gamma \leq 1.28057976 $ & $\gamma = 0.57281606 $ & $c_2 = 2 \gamma = 1.14563212
$ \\ \hline
\end{tabular}
\end{center}
\end{table}

\section{ESDIRK Integration Methods}
\label{sec:ESDIRKMethods}

In this section, we apply the order conditions
(\ref{eq:OrderConditionsRK}) and the conditions for A- and
L-stability to derive stiffly accurate ESDIRK methods of various
order. In addition, we equip these methods with continuous extensions.

\subsection{ESDIRK12}

The stability function of a two stage ESDIRK method is
\begin{equation}
    R(z) = \frac{1 + b_1 z}{1- \gamma z}
\end{equation}
L-stability requires the order of the numerator to be less than
the order of the denominator. Hence, L-stability gives the
requirement
\begin{equation}
    b_1 = 0 \qquad \qquad \text{(and $\gamma \neq 0$)}
\end{equation}
The order and consistency conditions for the two-stage ESDIRK
method becomes
\begin{subequations}
\begin{alignat}{3}
    &\text{Consistency / Order 1:} \qquad && b_1 + \gamma = 1 \\
    &\text{Order 2:} \qquad && \gamma = \frac{1}{2}
\end{alignat}
\end{subequations}
It is apparent that the only second-order method is $b_1 = \gamma
= \frac{1}{2}$ (the Trapez method). This method is A-stable, but
not L-stable ($b_1 \neq 0$). Hence, the maximum order of a
two-stage L-stable stiffly accurate ESDIRK method is 1. This
method has the coefficients $b_1 = 0$ and $\gamma = 1$, i.e. it is
the implicit Euler method. This method is also A-stable. The
second order method embedded in this implicit Euler method must
satisfy the conditions
\begin{subequations}
\begin{alignat}{3}
    & \text{Order 1:} \qquad && \hat{b}_1 +
    \hat{b}_2 = 1\\
    & \text{Order 2:} \qquad && \hat{b}_2 = \frac{1}{2}
\end{alignat}
\end{subequations}
The embedded method is uniquely determined as $\hat{b}_1 =
\hat{b}_2 = \frac{1}{2}$, i.e. as trapez quadrature. The embedded
method has the stability function
\begin{equation}
    \hat{R}(z) = \frac{1 - \frac{1}{2} z^2}{1 - z} \qquad |
    \hat{R}(\infty) | = \infty
\end{equation}
which is neither A- nor L-stable. However, this is of little
concern since it is the output of the basic integration method
that is used for the next step. The embedded method is merely used
to estimate the error provided within a single step.

Consequently, an A- and L-stable stiffly accurate ESDIRK method
with two-stages consists of the implicit Euler method as the basic
integrator and trapez quadrature for estimation of the error. The
ESDIRK12 method may be summarized by the Butcher tableau
\begin{equation}
\renewcommand{\arraystretch}{1.2}
    \begin{array}{l|ll}
        0 & 0 & \\
        1 & b_{1} & \gamma \\ \hline
          & b_{1} & \gamma \\
          & \hat{b}_1 & \hat{b}_2 \\ \hline
          & d_1 & d_2
    \end{array}
\qquad = \qquad
    \begin{array}{l|ll}
        0 & 0 & \\
        1 & 0 & 1 \\ \hline
          & 0 & 1 \\
          & \frac{1}{2} & \frac{1}{2} \\ \hline
          & - \frac{1}{2} & \frac{1}{2}
    \end{array}
\end{equation}
The continuous extension of ESDIRK12 can be of order 1 or 2,
respectively. In the case of order 1, the continuous extension
must satisfy the order 1 conditions. The additional degrees of
freedom is used to impose the conditions $b_i(\theta=1) = b_i$.
These conditions uniquely determines the order 1 continuous
extension as
\begin{equation}
    \begin{bmatrix} \bar{b}_1(\theta) \\ \bar{b}_2(\theta) \end{bmatrix}
    = \begin{bmatrix} \bar{b}_{11} \\ \bar{b}_{21} \end{bmatrix}
    \theta  \qquad \begin{bmatrix} \bar{b}_{11} \\ \bar{b}_{21}
    \end{bmatrix} = \begin{bmatrix} 0 \\ 1 \end{bmatrix}
\end{equation}
The order 2 continuous extension of ESDIRK12 is uniquely
determined by the conditions for order 1 and 2:
\begin{equation}
\renewcommand{\arraystretch}{1.2}
    \begin{bmatrix}
        \bar{b}_1(\theta) \\ \bar{b}_2(\theta)
    \end{bmatrix}
    =
    \begin{bmatrix}
        \bar{b}_{11} \\ \bar{b}_{21}
    \end{bmatrix}
    \theta
    + \begin{bmatrix} \bar{b}_{12} \\ \bar{b}_{22}
    \end{bmatrix} \theta^2 \qquad
    \begin{bmatrix}
        \bar{b}_{11} & \bar{b}_{12} \\
        \bar{b}_{21} & \bar{b}_{22}
    \end{bmatrix}
    = \begin{bmatrix} 1 & \frac{-1}{2} \\ 0 & \frac{1}{2}
    \end{bmatrix}
\end{equation}
It should be noted that $b_i(\theta = 1) = \hat{b}_i$ for $i=1,2$
in the case of the order 2 extension.

\subsection{ESDIRK23}
The stability function of the 3-stage stiffly accurate ESDIRK
integration scheme is
\begin{equation}
    R(z) = \frac{1 + \left( b_1 + b_2 - \gamma \right) z + \left( a_{21} b_2 - b_1 \gamma \right) z^2}{(1 - \gamma z)^2}
\end{equation}
To have L-stability, the numerator order must be less than the
denominator order in the stability function, i.e.
\begin{equation}
\label{eq:ESDIRK23:LstabilityCondition}
    a_{21} b_2 - b_1 \gamma = 0
\end{equation}
The consistency requirements for the ESDIRK23 scheme are
\begin{subequations}
\begin{alignat}{3}
    c_2 &= a_{21} + \gamma \label{eq:ESDIRK23:ConsistencyA}\\
    1   &= b_{1} + b_{2} + \gamma \label{eq:ESDIRK23:ConsistencyB}
\end{alignat}
\end{subequations}
and the order conditions of ESDIRK23 are
\begin{subequations}
\label{eq:ESDIRK23:Order}
\begin{alignat}{3}
    & \text{Order 1:} \qquad && b_1 + b_2 + \gamma = 1
    \label{eq:ESDIRK23:Order1}
\\  & \text{Order 2:} \qquad && b_2 c_2 + \gamma = \frac{1}{2} \label{eq:ESDIRK23:Order2}
\\  & \text{Order 3:} \qquad && b_2 c_2^2 + \gamma = \frac{1}{3} \label{eq:ESDIRK23:Order3a}
\\  &                        && 2 b_2 c_2 \gamma + \gamma^2 =
\frac{1}{6} \label{eq:ESDIRK23:Order3b}
\end{alignat}
\end{subequations}
The consistency condition (\ref{eq:ESDIRK23:ConsistencyB}) and the
order condition (\ref{eq:ESDIRK23:Order1}) are identical. Hence,
the consistency and order conditions for order 3 provide a system
of 5 nonlinear equations, (\ref{eq:ESDIRK23:ConsistencyA}) and
\ref{eq:ESDIRK23:Order}), in 5 unknown variables $\set{c_2,
a_{21}, \gamma, b_1, b_2}$. This system has two solutions
corresponding to $\gamma = \frac{3 \pm \sqrt{3}}{6}$.
None of them are L-stable.

Instead we aim at constructing a method of order 2. In addition we
require that stage 2 has order 2, i.e. that the conditions
\begin{subequations}
\label{eq:ESDIRK23:Stage2StageOrder2Conditions}
\begin{alignat}{3}
    & a_{21} + \gamma = c_2 \label{eq:ESDIRK23:StageOrder2A}\\
    &  \gamma c_2 = \frac{1}{2} c_2^2 \label{eq:ESDIRK23:StageOrder2B}
\end{alignat}
\end{subequations}
are satisfied. These conditions are equivalent to $c_2 = 2 \gamma$
and $a_{21} = \gamma$. Note that
(\ref{eq:ESDIRK23:Stage2StageOrder2Conditions}) and the order
conditions (\ref{eq:ESDIRK23:Order1}-\ref{eq:ESDIRK23:Order2})
automatically provide consistency. The L-stability condition,
(\ref{eq:ESDIRK23:LstabilityCondition}), the order conditions
(\ref{eq:ESDIRK23:Order1}-\ref{eq:ESDIRK23:Order2}), and the
conditions for stage order 2 of stage 2,
(\ref{eq:ESDIRK23:Stage2StageOrder2Conditions}), constitute 5
nonlinear equations in 5 unknown variables $\set{c_2, a_{21},
\gamma, b_1, b_2}$. This system has two solutions, corresponding
to $\gamma = \frac{2 \pm \sqrt{2}}{2}$. Both of them are A-stable.
However, only the solution corresponding to $\gamma = \frac{2 -
\sqrt{2}}{2}$ has $0 < c_2 < 1$. The additional requirement $0 <
c_2 < 1$ thus provides the unique solution: $\gamma = a_{21} =
\frac{2-\sqrt{2}}{2} \approx 0.2929$, $c_2 = 2 \gamma = 2 -
\sqrt{2} \approx 0.5858$, $b_1 = b_2 = \frac{1-\gamma}{2} =
\frac{\sqrt{2}}{4} \approx 0.3536$.

Having a stiffly accurate, L-stable 3-stage ESDIRK integration
scheme of order 2, an embedded method of order 3 may be determined
as the solution of the following conditions:
\begin{subequations}
\begin{alignat}{3}
    & \text{Order 1:} \qquad && \hat{b}_1 + \hat{b}_2 + \hat{b}_3
    = 1 \\
    & \text{Order 2:} \qquad && \hat{b}_2 c_2 + \hat{b}_3 =
    \frac{1}{2} \\
    & \text{Order 3:} \qquad && \hat{b}_2 c_2^2 + \hat{b}_3 =
    \frac{1}{3}     \label{eq:ESDIRK23:EmbeddedOrderConditions3a} \\
    & && \hat{b}_2 (c_2 \gamma) + \hat{b}_3 (c_2 b_2 + \gamma) =
    \frac{1}{6}
    \label{eq:ESDIRK23:EmbeddedOrderConditions3b}
\end{alignat}
\end{subequations}
These order conditions constitute 4 linear equations with 3
unknown variables, $\set{\hat{b}_1, \hat{b}_2, \hat{b}_3}$.
However, it turns out that
(\ref{eq:ESDIRK23:EmbeddedOrderConditions3b}) is linearly
dependent of (\ref{eq:ESDIRK23:EmbeddedOrderConditions3a}) as
$c_2^2 = 4 \gamma^2$, $c_2\gamma = 2 \gamma^2 = \frac{1}{2}
c_2^2$, and $c_2 b_2 + \gamma = \frac{1}{2}$ by condition
(\ref{eq:ESDIRK23:Order2}). Consequently, the variables,
$\set{\hat{b}_1, \hat{b}_2, \hat{b}_3}$, and $\set{d_i = b_i -
\hat{b}_i}_{i=1}^{3}$ can be uniquely determined as:
\begin{subequations}
\begin{alignat}{3}
    \hat{b}_1 &= \frac{6 \gamma -1}{12 \gamma} \approx 0.2155
     && d_1 =  \frac{1-6 \gamma^2}{12
    \gamma} \approx 0.1381
\\    \hat{b}_2 &= \frac{1}{12 \gamma (1 - 2 \gamma)} \approx
0.6869 \qquad  && d_2 = \frac{6\gamma(1-2\gamma)(1-\gamma)-1}{12
\gamma (1-2\gamma)} \approx -0.3333
\\  \hat{b}_3 &= \frac{1 - 3
\gamma}{3 (1-2\gamma)} \approx 0.0976  && d_3 = \frac{6 \gamma
(1-\gamma)-1}{3 (1-2\gamma)} \approx 0.1953
\end{alignat}
\end{subequations}
The transfer function of the embedded method is
\begin{equation}
    \hat{R}(z) = \frac{ \frac{-10 + 7 \sqrt{2}}{6(\sqrt{2}-1)} z^3
    + \frac{3-2\sqrt{2}}{\sqrt{2}-1} z + 1}{\left(1-
    \frac{2-\sqrt{2}}{2} z \right)^2} \qquad | \hat{R}(\infty) | =
    \infty
\end{equation}
which is neither A- nor L-stable. This ESDIRK method, called
ESDIRK23, may be summarized by the Butcher tableau:
\begin{equation*}
\renewcommand{\arraystretch}{1.2}
\begin{array}{l|lll}
    0 & 0 & & \\
    c_2 & a_{21} & \gamma & \\
    1   & b_1 & b_2 & \gamma \\ \hline
        & b_1 & b_2 & \gamma \\
        & \hat{b}_1 & \hat{b}_2 & \hat{b}_3 \\ \hline
        & d_1 & d_2 & d_3
\end{array}
\,\, = \,\,
\begin{array}{c|ccc}
    0 & 0 & & \\
    2\gamma & \gamma & \gamma & \\
    1   & \frac{1-\gamma}{2} & \frac{1-\gamma}{2} & \gamma \\ \hline
        & \frac{1-\gamma}{2} & \frac{1-\gamma}{2} & \gamma \\
        &  \frac{6 \gamma -1}{12 \gamma}  & \frac{1}{12\gamma(1-2\gamma)} & \frac{1-3\gamma}{3(1-2\gamma)} \\ \hline
        & \frac{1-6 \gamma^2}{12
    \gamma} & \frac{6\gamma(1-2\gamma)(1-\gamma)-1}{12
\gamma (1-2\gamma)}  & \frac{6 \gamma (1-\gamma)-1}{3 (1-2\gamma)}
\end{array}
\qquad \gamma = \frac{2-\sqrt{2}}{2}
\end{equation*}
The continuous extension of order 2 with $\bar{x}(t_n+h) = X_3 =
x_{n+1}$ cannot be determined uniquely but has one-degree of
freedom left. If we use the spare degree of freedom to satisfy the
additional requirement $\bar{x}(t_n+c_2h) = X_2$ we get the
following unique 2nd order continuous extension
\begin{equation}
\renewcommand{\arraystretch}{1.2}
    \begin{bmatrix}
        \bar{b}_1(\theta) \\ \bar{b}_2(\theta) \\ \bar{b}_3(\theta)
    \end{bmatrix}
    =
    \begin{bmatrix}
        \bar{b}_{11} \\ \bar{b}_{21} \\ \bar{b}_{31}
    \end{bmatrix} \theta
    +
    \begin{bmatrix}
        \bar{b}_{12} \\ \bar{b}_{22} \\ \bar{b}_{32}
    \end{bmatrix} \theta^2
\qquad
    \begin{bmatrix}
        \bar{b}_{11} & \bar{b}_{12} \\
        \bar{b}_{21} & \bar{b}_{22} \\
        \bar{b}_{31} & \bar{b}_{32}
    \end{bmatrix}
    =
    \begin{bmatrix}
        \frac{\sqrt{2}}{2} & \frac{-\sqrt{2}}{4}\\
        \frac{\sqrt{2}}{2} & \frac{-\sqrt{2}}{4} \\
        1-\sqrt{2} & \frac{\sqrt{2}}{2}
    \end{bmatrix}
\end{equation}
The continuous extension of order 3 is given uniquely and is
%\begin{subequations}
\begin{equation}
    \begin{bmatrix}
        \bar{b}_1(\theta) \\ \bar{b}_2(\theta) \\ \bar{b}_3(\theta)
    \end{bmatrix}
    =
    \begin{bmatrix}
        \bar{b}_{11} \\ \bar{b}_{21} \\ \bar{b}_{31}
    \end{bmatrix} \theta
    +
    \begin{bmatrix}
        \bar{b}_{12} \\ \bar{b}_{22} \\ \bar{b}_{32}
    \end{bmatrix} \theta^2
    +
    \begin{bmatrix}
        \bar{b}_{13} \\ \bar{b}_{23} \\ \bar{b}_{33}
    \end{bmatrix} \theta^3
\end{equation}
in which
\begin{equation*}
\qquad
\begin{bmatrix}
    \bar{b}_{11} & \bar{b}_{12} & \bar{b}_{13} \\
    \bar{b}_{21} & \bar{b}_{22} & \bar{b}_{23} \\
    \bar{b}_{31} & \bar{b}_{32} & \bar{b}_{33}
\end{bmatrix}
=
\begin{bmatrix}
    1 & -1.35355339059327 &  0.569035593728849 \\
    0 &  2.06066017177982 &  -1.37377344785321 \\
    0 &  -0.707106781186547 & 0.804737854124365
\end{bmatrix}
\end{equation*}
%\end{subequations}
It satisfies the 3rd order conditions and $\bar{x}(t_n+h) =
\hat{x}_{n+1}$. It is not possible to construct a 3rd order
continuous extension satisfying $\bar{x}(t_n+h) = x_{n+1} = X_3$.

\subsection{ESDIRK34}
We have developed the methods ESDIRK12 and ESDIRK23 in a quite
detailed way. The same can be done in the development of ESDIRK34.
However, we will develop the method in a direct way applying the
results of Table \ref{table:StabilityRegionsESDIRK}. There exists
no A- and L-stable stiffly accurate ESDIRK method with 4 stages of
order 4. According to Table \ref{table:StabilityRegionsESDIRK},
the diagonal coefficient $\gamma = 0.43586652$ of the A- and
L-stable stiffly accurate ESDIRK method is unique. In the
following we will determine the coefficients of ESDIRK34 such that
it is A- and L-stable, stiffly accurate with order 3 of the
advancing method and order 4 of the embedded method. Continuous
extensions to this method will be developed as well.

The stability function of a 4 stage stiffly accurate ESDIRK
integration method is
\begin{equation}
\begin{split}
    R(z) &= \frac{\left[a_{21} a_{32} b_3 - (a_{21} b_2 + a_{31} b_3 ) \gamma + b_1 \gamma^2
    \right]z^3}{(1-\gamma z)^3} \\ & \quad + \frac{
    \left[ (a_{21}b_2 + a_{31} b_3 + a_{32} b_3) - (2b_1+b_2+b_3)\gamma + \gamma^2 \right] z^2 + \left[b_1 + b_2 + b_3 - 2\gamma \right] z + 1}{(1-\gamma z)^3}
\end{split}
\end{equation}
which implies that a requirement for L-stability is
\begin{equation}
\label{eq:ESDIRK34Lstability}
    a_{21} a_{32} b_3 - (a_{21} b_{2} + a_{31} b_{3}) \gamma + b_1
    \gamma^2 = 0
\end{equation}
The consistency conditions (\ref{eq:ConsistencyConditionRK}) are
\begin{subequations}
\label{eq:ESDIRK34consistency}
\begin{alignat}{3}
    c_2 &= a_{21} + \gamma \\
    c_3 &= a_{31} + a_{32} + \gamma \\
    1 &= b_1 + b_2 + b_3 + \gamma
\end{alignat}
\end{subequations}
and the conditions for order 3 of the advancing method are
\begin{subequations}
\label{eq:ESDIRK34AdvancingOrder3}
\begin{alignat}{3}
    & \text{Order 1:} \qquad && b_1 + b_2 + b_3 + \gamma = 1 \\
    & \text{Order 2:} \qquad && b_2 c_2 + b_3 c_3 + \gamma =
    \frac{1}{2} \\
    & \text{Order 3:} \qquad && b_2 c_2^2 + b_3 c_3^2+ \gamma =
    \frac{1}{3} \\
    & && b_3 a_{32} c_2 + 2(b_2 c_2 + b_3 c_3)\gamma + \gamma^2 = \frac{1}{6}
\end{alignat}
\end{subequations}
Stage order 2 for the stages 2 and 3, i.e. $\tilde{C}(2)$
(\ref{eq:ModifiedStageOrder2}), gives the additional relations
\begin{subequations}
\label{eq:ESDIRK34ModifiedStageOrder2}
\begin{alignat}{3}
    & \gamma c_2 = \frac{1}{2} c_2^2
    \\& a_{32} c_2 + \gamma c_3 = \frac{1}{2} c_3^2
\end{alignat}
\end{subequations}
Along with the requirement of A-stability, i.e. $\frac{1}{3} \leq
\gamma \leq 1.06857902$, the advancing method is uniquely
determined by the conditions
(\ref{eq:ESDIRK34Lstability})-(\ref{eq:ESDIRK34ModifiedStageOrder2})
\cite{Alexander:2003}.

When the coefficients of the advancing method have been
determined, the conditions for order 4 gives the following linear
relations that $\hat{b}$ must satisfy
\begin{subequations}
\begin{alignat}{3}
    & \text{Order 1:} \qquad && \hat{b}_1 + \hat{b}_2 + \hat{b}_3
    + \hat{b}_4 = 1
\\  & \text{Order 2:} \qquad &&  \hat{b}_2 c_2 + \hat{b}_3 c_3 +
\hat{b}_4 = \frac{1}{2}
\\  & \text{Order 3:} \qquad && \hat{b}_2 c_2^2 + \hat{b}_3 c_3^2
+ \hat{b}_4 = \frac{1}{3}
\\ & && \hat{b}_2 \left( c_2 \gamma \right) + \hat{b}_3
\left(a_{32} c_2 + \gamma c_3 \right) + \hat{b}_4 \left( b_2 c_2 +
b_3 c_3 + \gamma \right) = \frac{1}{6}
\\  & \text{Order 4:} \qquad && \hat{b}_2 c_2^3 + \hat{b}_3 c_3^3
+ \hat{b}_4 = \frac{1}{4}
\\ & && \hat{b}_2 \gamma c_2^2 + \hat{b}_3 c_3 (a_{32} c_2 +
\gamma c_3)+ \hat{b}_4 (b_2 c_2 + b_3 c_3 + \gamma) = \frac{1}{8}
\\ & && \hat{b}_2 \gamma c_2^2 + \hat{b}_3 (a_{32} c_2^2 + \gamma
c_3^2) + \hat{b}_4 (b_2 c_2^2 + b_3 c_3^2 + \gamma) = \frac{1}{12}
\\ & &&
 \hat{b}_2 \gamma^2 c_2 + \hat{b}_3 (2 a_{32} \gamma c_2 +
\gamma^2 c_3) \nonumber \\ & && \qquad \qquad  + \hat{b}_4 \left(
(2 b_2 \gamma + b_3 a_{32}) c_2 + 2 b_3 \gamma c_3 + \gamma^2
\right) = \frac{1}{24}
\end{alignat}
\end{subequations}
or in more compact notation
\begin{equation}
\renewcommand{\arraystretch}{1.2}
    \left[
    \begin{array}{c}
        e' \\ \hline
        (Ce)' \\ \hline
        (C^2 e)' \\
        (ACe)' \\ \hline
        (C^3 e)' \\
        (CACe)' \\
        (AC^2 e)' \\
        (A^2 C e)'
    \end{array}
    \right]
    \begin{bmatrix} \hat{b}_1 \\ \hat{b}_2 \\ \hat{b}_3 \\
    \hat{b}_4 \end{bmatrix}
    = \left[ \begin{array}{c} 1 \\ \hline \frac{1}{2} \\ \hline \frac{1}{3} \\
    \frac{1}{6} \\ \hline \frac{1}{4} \\ \frac{1}{8} \\ \frac{1}{12} \\
    \frac{1}{24} \end{array} \right]
\end{equation}
The solution, $\hat{b} =
\begin{bmatrix} \hat{b}_1 & \hat{b}_2 & \hat{b}_3 &
\hat{b}_4\end{bmatrix}'$, to  this over-determined linear system
exists and is unique. The coefficients of the error estimator is
determined as $d = b - \hat{b}$. The embedded method has the
stability function
\begin{equation}
    \hat{R}(z) = \frac{1-0.3076 z - 0.2377 z^2 + 0.2590 z^4}{\left(1-0.4359 z\right)^3} \qquad |\hat{R}(\infty)| = \infty
\end{equation}
which is neither A- nor L-stable.

The developed 4-stage, stiffly accurate, A- and L-stable ESDIRK
method of third order with an embedded method of order 4 is called
ESDIRK34. It is summarized by the Butcher tableau
\begin{equation*}
\begin{array}{l|l}
    c & A \\ \hline & b' \\ & \hat{b}' \\ \hline & d'
\end{array} \quad = \quad
\begin{array}{l|llll}
    0   & 0 & & & \\
    c_2 & a_{21} & \gamma & & \\
    c_3 & a_{31} & a_{32} & \gamma & \\
    1   & b_1 & b_2 & b_3 & \gamma \\ \hline
        & b_1 & b_2 & b_3 & \gamma \\
        & \hat{b}_1 & \hat{b}_2 & \hat{b}_3 & \hat{b}_4 \\ \hline
        & d_1 & d_2 & d_3 & d_4
\end{array}
\end{equation*}
with the coefficients listed in Table
\ref{table:ESDIRK34coefficients}. In particular the location of
the quadrature points should be noted, i.e. $0 = c_1 < c_3 < c_2 <
c_4 = 1$ which shows that $c_2 > c_3$.
\begin{table}[tb]
\caption{Coefficients for ESDIRK34}.
\label{table:ESDIRK34coefficients} \footnotesize
\begin{center}
\begin{tabular}{lccc} \hline
    $i$ & $b_i$ & $\hat{b}_i$ & $d_i$ \\ \hline
    1 & 0.10239940061991099768 & 0.15702489786032493710 & -0.05462549724041393942\\
    2 & -0.3768784522555561061 & 0.11733044137043884870 & -0.49420889362599495480\\
    3 & 0.83861253012718610911 & 0.61667803039212146434 &  0.22193449973506464477\\
    4 & 0.43586652150845899942 & 0.10896663037711474985 &  0.32689989113134424957\\ \hline
      & $a_{21}$ & $a_{31}$ & $a_{32}$  \\ \hline
      & 0.43586652150845899942 & 0.14073777472470619619 &
      -0.1083655513813208000 \\ \hline
      & $\gamma$ & $c_2$ & $c_3$ \\ \hline
      & 0.43586652150845899942 & 0.87173304301691799883 &
      0.46823874485184439565 \\ \hline
\end{tabular}
\end{center}
\end{table}

Using the procedure introduced in Section \ref{sec:ContinuousExtension},
construction of different continuous extensions has been attempted.
There exists no 2nd order continuous extension
\begin{equation}
    \bar{b}(\theta) = \bar{b}_1 \theta + \bar{b}_2 \theta^2 = \bar{B} \begin{bmatrix} \theta \\ \theta^2
    \end{bmatrix} \qquad \bar{B} = \begin{bmatrix} \bar{b}_1 &
    \bar{b}_2 \end{bmatrix}
\end{equation}
satisfying $\bar{x}(t_n+c_ih) = X_i$ for $i=2,3,4$. There exists a
unique 2nd order continuous extension
\begin{equation*}
    \bar{B}_{24} = \begin{bmatrix} \bar{b}_1 & \bar{b}_2 \end{bmatrix}
    = \begin{bmatrix}
          3.20218915732655 & -3.09978975670664 \\
          6.45947654423207 & -6.83635499648762 \\
         -5.69941214787150  & 6.53802467799868 \\
         -2.96225355368712 &  3.39812007519558
     \end{bmatrix}
\end{equation*}
satisfying $\bar{x}(t_n+c_ih) = X_i$ for $i=2,4$ and another
unique 2nd order continuous extension
\begin{equation*}
    \bar{B}_{34} = \begin{bmatrix} \bar{b}_1 & \bar{b}_2 \end{bmatrix}
    = \begin{bmatrix}
         0.47506477777383     &   -0.372665377153919 \\
        -0.103360609602923    &    -0.273517842652633 \\
          1.01209512329345    &    -0.173482593166265 \\
        -0.383799291464359    &     0.819665812972817
    \end{bmatrix}
\end{equation*}
satisfying $\bar{x}(t_n+c_ih) = X_i$ for $i=3,4$. Obviously, no
3rd order continuous extension
\begin{equation}
    \bar{b}(\theta) = \bar{b}_1 \theta + \bar{b}_2 \theta^2 +
    \bar{b}_3 \theta^3 = \bar{B} \begin{bmatrix} \theta \\
    \theta^2 \\ \theta^3 \end{bmatrix} \quad \bar{B} =
    \begin{bmatrix} \bar{b}_1 & \bar{b}_2 & \bar{b}_3
    \end{bmatrix}
\end{equation}
satisfies $\bar{x}(t_n+c_ih) = X_i$ for $i=2,3,4$ as no 2nd order
continuous extension does so. Furthermore, it is not possible to
construct continuous extensions of order 3 satisfying
$\bar{x}(t_n+c_ih)=X_i$ for either $i=2,4$ or $i=3,4$. In
contrast, the 3rd order continuous extension satisfying
$\bar{x}(t_n+h) = X_4 = x_{n+1}$ is not unique. The coefficient
matrix of one such continuous extension is (the minimum norm solution obtained using SVD)
\begin{equation*}
     \begin{bmatrix} \bar{b}_1 & \bar{b}_2 & \bar{b}_3
    \end{bmatrix}
    = \begin{bmatrix}
         0.969611875176691 &         -1.53835725968354 &        0.671144785126761 \\
        -0.274928052044991 &        0.266658367468879  &      -0.368608767679444 \\
         0.123462002567514 &         1.88835458133267  &         -1.173204053773 \\
         0.181854174300786 &        -0.61665568911801  &       0.870668036325683
    \end{bmatrix}
\end{equation*}
and
\begin{equation*}
    \begin{bmatrix} \bar{b}_1 & \bar{b}_2 & \bar{b}_3 \end{bmatrix} =
    \begin{bmatrix}
               0.927166003679448 &         -1.64140141649749 &         0.816634813437953 \\
        -0.658945191501327   &     -0.665604793624494       &  0.947671532870265 \\
          0.29591266631342   &       2.30700621012198       &  -1.76430634630822 \\
          0.43586652150846   & 0                       &  0
    \end{bmatrix}
\end{equation*}
for a continuous extension that in addition has minimum curvature of $\bar{b}_4(\theta)$.
The minimum norm 3rd order continuous extension satisfying $\bar{x}(t_n+h) = X_4$ and
$\dot{\bar{x}}(t_n+h) = \dot{X}_4$ has the coefficients
\begin{equation*}
    \begin{bmatrix} \bar{b}_1 & \bar{b}_2 & \bar{b}_3 \end{bmatrix} =
    \begin{bmatrix}
   0.92277773077164 &  -1.53835725968353 &   0.71797892953181 \\
  -0.69864686211777 &   0.26665836746888 &   0.05511004239334 \\
   0.31374150452444 &   1.88835458133266 &  -1.36348355572992 \\
   0.46212762682169 &  -0.61665568911801 &   0.59039458380477
\end{bmatrix}
\end{equation*}
There exists no continuous extension satisfying the conditions for
order 4.

\section{ESDIRK Methods Family due to Kv{\ae}rn{\o}}
\label{sec:ESDIRKKvaernoe} Kv{\ae}rn{\o}~\cite{Kvaernoe:2004}
considers a class of ESDIRK methods in which both the advancing
method and the embedded method are stiffly accurate and A-stable.
The advancing method is also L-stable while $|\hat{R}(\infty)|$
for the embedded method is minimized. These extra properties come
at the expense of more implicit stages to attain methods of a
given order. The structure of the Butcher tableaus for the methods
developed by Kv{\ae}rn{\o} is illustrated by the Butcher tableau
for a 5-stage method
\begin{equation}
\label{eq:ESDIRKKvaernoe}
    \begin{array}{c|c}
        c & A \\ \hline
          & b' \\
          & \hat{b}' \\ \hline & d'
    \end{array}
    \quad = \qquad
    \begin{array}{c|ccccc}
        0 & 0 & 0 & 0 & 0 & 0 \\
        c_2 & a_{21} & \gamma & 0 & 0 & 0 \\
        c_3 & a_{31} & a_{32} & \gamma & 0 & 0 \\
        1   & b_1 & b_2 & b_3 & \gamma & 0 \\
        1   & \hat{b}_1 & \hat{b}_2 & \hat{b}_3 & \hat{b}_4 &
        \gamma \\ \hline
         & b_1 & b_2 & b_3 & \gamma & 0 \\
         & \hat{b}_1 & \hat{b}_2 & \hat{b}_3 & \hat{b}_4 & \gamma
         \\ \hline
         & d_1 & d_2 & d_3 & d_4 & d_5
    \end{array}
\end{equation}

\subsection{ESDIRK 3/2 with 4 stages}

The Butcher tableau for Kv{\ae}rn{\o}'s ESDIRK method of order 3/2
with 4 stages is
\begin{equation*}
\renewcommand{\arraystretch}{1.5}
    \begin{array}{c|cccc}
        0 & 0 & 0 & 0 & 0 \\
        c_2 & a_{21} & \gamma & 0 & 0 \\
        1   & \hat{b}_1 & \hat{b}_2 & \gamma & 0 \\
        1  & b_1 & b_2 & b_3 & \gamma \\ \hline
           & b_1 & b_2 & b_3 & \gamma \\
           & \hat{b}_1 & \hat{b}_2 & \gamma & 0 \\ \hline
           & d_1 & d_2 & d_3 & d_4
    \end{array}
    \quad = \quad
    \begin{array}{c|cccc}
        0 & 0 & 0 & 0 & 0 \\
        2\gamma & \gamma & \gamma & 0 & 0 \\
        1 & \frac{-4\gamma^2 +6\gamma-1}{4\gamma} &
        \frac{-2\gamma+1}{4\gamma} & \gamma & 0 \\
        1 & \frac{6\gamma-1}{12\gamma} &
        \frac{-1}{12\gamma(2\gamma-1)} &
        \frac{-6\gamma^2+6\gamma-1}{3(2\gamma-1)} & \gamma \\ \hline
        & \frac{6\gamma-1}{12\gamma} &
        \frac{-1}{12\gamma(2\gamma-1)} &
        \frac{-6\gamma^2+6\gamma-1}{3(2\gamma-1)} & \gamma  \\
        & \frac{-4\gamma^2 +6\gamma-1}{4\gamma} &
        \frac{-2\gamma+1}{4\gamma} & \gamma & 0 \\ \hline
        & \frac{6 \gamma(\gamma-1) +1}{6 \gamma} &
        \frac{3(2\gamma-1)^2-1}{12\gamma(2\gamma-1)} &
        \frac{-12\gamma^2 + 9 \gamma -1}{3(2\gamma-1)} & \gamma
    \end{array}
\end{equation*}
In the case $x_{n+1} = X_4$, $\gamma = 0.4358665215$ and
$|\hat{R}(\infty)|=0.9569$. This method is called ESDIRK32a. No
continuous extension of order 3 being identical with the advancing
method in the end-point and in the internal point exists. The
minimum norm continuous extension of order 3 satisfying
$\bar{x}(t_n+h) = X_4$ and $\dot{\bar{x}}(t_n+h) = \dot{X}_4$ for
ESDIRK32a has the coefficients
\begin{equation*}
    \begin{bmatrix}
        \bar{b}_1 & \bar{b}_2 & \bar{b}_3
    \end{bmatrix}
    =
    \begin{bmatrix}
    1.00000000000000 & -1.07357009006975 &   0.38238006004650 \\
    0.00000000000000 &  4.47169016526534 & -2.98112677684356  \\
   -0.86407093427697 & -1.97757777116702 &  1.60640882553700  \\
    0.86407093427697 & -1.42054230402855 &  0.99233789126005  \\
    \end{bmatrix}
\end{equation*}

In the case $x_{n+1} = X_3$, $\gamma = \frac{2-\sqrt{2}}{2}$ and
$|\hat{R}(\infty)| = 1.609 > 1$. This method is called ESDIRK32b.
The unique continuous extension of order 2 for ESDIRK32b
satisfying $\bar{x}(t_n+c_2h) = X_2$, $\bar{x}(t_n+h) = X_3$, and
$\dot{\bar{x}}(t_n+h) = \dot{X}_3$ has the coefficients
\begin{equation*}
\renewcommand{\arraystretch}{1.2}
    \begin{bmatrix}
        \bar{b}_1 & \bar{b}_2
    \end{bmatrix}
    =
    \begin{bmatrix}
        \frac{\sqrt{2}}{2} & - \frac{\sqrt{2}}{4} \\
        \frac{\sqrt{2}}{2} & - \frac{\sqrt{2}}{4} \\
        1-\sqrt{2} & \frac{\sqrt{2}}{2} \\
        0 & 0
    \end{bmatrix}
\end{equation*}

\subsection{ESDIRK 4/3 with 5 stages}
This method is represented by Butcher tableau
(\ref{eq:ESDIRKKvaernoe}) and has coefficients provided in
\cite{Kvaernoe:2004}.

In the case $x_{n+1} = X_5$, $\gamma = 0.5728160625$ and
$|\hat{R}(\infty)| = 0.5525$. This method is called ESDIRK43a
However, $c_2 = 2 \gamma > 1$. Therefore, we disregard this method
as it is not useful as a general purpose ESDIRK integration
algorithm applicable to discrete event systems.

In the case $x_{n+1} = X_4$, $\gamma = 0.4358665215$ and
$|\hat{R}(\infty)| = 0.7175$. This method is called ESDIRK43b. The
matrices in the Butcher tableau of ESDIRK43b are {\scriptsize
\begin{align*}
&A = \\&\begin{bmatrix}
                  0 &                 0 &                 0 &                 0  &                0  \\
   0.43586652150846 &  0.43586652150846 &                 0 &                 0  &                0  \\
   0.14073777472471 & -0.10836555138132 &  0.43586652150846 &                 0  &                0  \\
   0.10239940061991 & -0.37687845225556 &  0.83861253012719 &  0.43586652150846  &                0  \\
   0.15702489786032 &  0.11733044137044 &  0.61667803039212 & -0.32689989113134  & 0.43586652150846
    \end{bmatrix}
\\
&b' =
    \begin{bmatrix}
   0.10239940061991 &  -0.37687845225556  & 0.83861253012719 &  0.43586652150846 &                 0
    \end{bmatrix}
\\
&\hat{b}' =
\begin{bmatrix}
   0.15702489786032 &   0.11733044137044 &  0.61667803039212 & -0.32689989113134 &  0.43586652150846
\end{bmatrix}
\\
&c' = \begin{bmatrix}
                 0 & 0.87173304301692 & 0.46823874485185 & 1 & 1
\end{bmatrix}
\end{align*}
} The corresponding minimum norm continuous extension of order 3
satisfying $\bar{x}(t_n+h) = X_4$ and $\dot{\bar{x}}(t_n+h) =
\dot{X}_4$ has the coefficients
\begin{equation*}
    \begin{bmatrix} \bar{b}_1 & \bar{b}_2 & \bar{b}_3
    \end{bmatrix}
    =
    \begin{bmatrix}
   0.91305667617487 &  -1.51891515049001 &  0.70825787493505 \\
  -0.78659538212849 &  0.44255540749030  & -0.03283847761737 \\
   0.35323656631463 &  1.80936445775230  & -1.32398849393974 \\
   0.30072875082513 & -0.29385793712489  &  0.42899570780821 \\
   0.21957338881385 & -0.43914677762771  &  0.21957338881385
    \end{bmatrix}
\end{equation*}
This method is called ESDIRK43b. There exists no continuous
extension of order 3 that in addition to the above conditions is
equal to the internal stage values of ESDIRK34, i.e. satisfies
$\bar{x}(t_n+ c_2 h) = X_2$ or/and $\bar{x}(t_n + c_3 h) = X_3$.
This non-existence observation holds even if the condition
$\dot{\bar{x}}(t_n+h) = \dot{X}_4$ is relaxed.

\subsection{ESDIRK 5/4 with 7 stages}
Kv{\ae}rn{\o} \cite{Kvaernoe:2004} provides two embedded ESDIRK
methods with 7 stages. They are of order 5 and 4, respectively.
For both methods, the advancing method is L-stable and stiffly
accurate, while the embedded method for error estimation is
A-stable and stiffly accurate.

We will not pay further consideration to these methods as linear
multi-step methods are usually preferable for high-accuracy
solutions.

\section{Other ESDIRK Methods}
\label{sec:OtherESDIRKMethods}

Williams {\em et.al}~\cite{Williams:Burrage:Cameron:Kerr:2002}
constructed an ESDIRK method of order 3. This method is
constructed such that it is applicable to index-2 differential
algebraic systems. This method is represented by the Butcher
tableau
\begin{equation*}
\renewcommand{\arraystretch}{1.2}
    \begin{array}{c|cccc}
        0 & \\
        c_2 & a_{21} & \gamma\\
        c_3 & a_{31} & a_{32} & \gamma\\
        1   & b_{1} & b_{2} & b_{3} & \gamma \\ \hline
            & b_{1} & b_{2} & b_{3} & \gamma \\
            & \hat{b}_{1} & \hat{b}_{2} & \hat{b}_{3} & \hat{b}_4
            \\ \hline
            & d_1 & d_2 & d_3 & d_4
    \end{array}
    \quad = \quad
    \begin{array}{c|cccc}
        0 \\
        1 & \frac{1}{2} & \frac{1}{2} \\
        \frac{3}{2} & \frac{5}{8} & \frac{3}{8} & \frac{1}{2} \\
        1 & \frac{7}{18} & \frac{1}{3} & -\frac{2}{9} & \frac{1}{2}
        \\ \hline
        & \frac{7}{18} & \frac{1}{3} & -\frac{2}{9} & \frac{1}{2}
        \\
        & \frac{1}{2} & \frac{1}{2} &  &  \\ \hline
        & -\frac{1}{9} &  -\frac{1}{6} & -\frac{2}{9} &
        \frac{1}{2}
    \end{array}
\end{equation*}
and we call it ESDIRK32c. However, this method is not suitable as
a general purpose method applicable to discrete-event systems as
$c_3 = \frac{3}{2} > 1$.

Butcher and Chen \cite{Butcher:Chen:2000} construct a 4th order A-
and L-stable method with stage order 2. The error estimator of
their method is {\em close} to 5th order. The method has 6 stages,
which is the minimum number of stages to have a 4th L-stable
ESDIRK method. We call this method ESDIRK45c. The Butcher-tableau
of this method is
\begin{equation*}
\renewcommand{\arraystretch}{1.2}
\begin{array}{l|cccccc}
    0 & 0 \\
    c_2 & a_{21} & \gamma \\
    c_3 & a_{31} & a_{32} & \gamma \\
    c_4 & a_{41} & a_{42} & a_{43} & \gamma \\
    c_5 & a_{51} & a_{52} & a_{53} & a_{54} & \gamma \\
    1 & b_1 & b_2 & b_3 & b_4 & b_5 & \gamma \\ \hline
      & b_1 & b_2 & b_3 & b_4 & b_5 & \gamma \\ \hline
      & d_1 & d_2 & d_3 & d_4 & d_5 & d_6
\end{array}
\quad = \quad
\begin{array}{l|cccccc}
    0 & 0 \\
    \frac{1}{2} & \frac{1}{4} & \frac{1}{4} \\
    \frac{1}{4} & \frac{1}{16} & \frac{-1}{16} & \frac{1}{4} \\
    \frac{1}{2} & \frac{-7}{36} & \frac{-4}{9} & \frac{8}{9} &
    \frac{1}{4} \\
    \frac{3}{4} & \frac{-5}{48} & \frac{-257}{768} & \frac{5}{6} &
    \frac{27}{256} & \frac{1}{4} \\
    1 & \frac{1}{4} & \frac{2}{3} & \frac{-1}{3} & \frac{1}{2} &
    \frac{-1}{3} & \frac{1}{4} \\ \hline
     & \frac{1}{4} & \frac{2}{3} & \frac{-1}{3} & \frac{1}{2} &
    \frac{-1}{3} & \frac{1}{4} \\ \hline
    & \frac{7}{90} & \frac{3}{20} & \frac{16}{45} & \frac{-1}{60}
    & \frac{16}{45} & \frac{7}{90}
\end{array}
\end{equation*}
However, since the embedded method of ESDIRK45c is not of order 5,
the behavior of error estimators and step size controllers are
uncertain. Due to such implementation considerations, we do not
give further consideration to ESDIRK45c.

\section{Conclusion}
\label{sec:conclusion}
%\subsection{Summary of ESDIRK Methods}

The properties of the ESDIRK methods discussed are summarized in
Table \ref{table:ESDIRKMethods}. The advancing method is an all
cases stiffly accurate as well as A- and L-stable. ESDIRK43a and
ESDIRK32c are not suitable for discrete-event systems as some of
the quadrature points are outside the interval of the current
step. ESDIRK45c is disregarded as the order of the embedded method
is uncertain. This yields unpredictable behavior of the step size
controller in an implementation of the method. ESDIRK54a and
ESDIRK54b are high order methods intended to obtain solutions of
high precision. Linear multi-step methods are usually regarded
most suitable for such integration tasks. The remaining ESDIRK
methods have been equipped with continuous extensions such that
they can be applied to discrete-event systems. They are suitable
to obtain low to medium accuracy solutions of stiff systems of
ordinary differential equations as well as systems of index-1
differential equations.

\begin{table}
\caption{Properties of the presented ESDIRK methods. All ESDIRK
integrators considered have an advancing method that is stiffly
accurate as well as A- and L-stable. s: Number of stages. $p$ and
$\hat{p}$: Order. A-S.: A-stablity. S. A.: Stiffly accurate.}
\label{table:ESDIRKMethods} \footnotesize
\begin{center}
\renewcommand{\arraystretch}{1.5}
\begin{tabular}{lccccccccccc}
\hline & &  & \multicolumn{4}{c}{Advancing Method} & &
\multicolumn{4}{c}{Embedded Method}  \\ \cline{4-7}\cline{9-12}
Method & s  & $\gamma$ & p & A-S. & $|R(\infty)|$ & S. A. &  &
$\hat{p}$ & A-S. & $|\hat{R}(\infty)|$ & {S. A.} \\ \hline
ESDIRK12 & 2 & 1 & 1 &  Yes & 0 & Yes  & &  2 & No & $\infty$ & No \\
ESDIRK23 & 3 & 0.2929  & 2 & Yes & 0 & Yes & &   3 & No & $\infty$ & No  \\
ESDIRK34 & 4 &  0.4359  & 3 & Yes & 0 & Yes  & &  4 & No & $\infty$ & No \\
\hline
ESDIRK32a & 4 & 0.4359 & 3  & Yes & 0 & Yes  &  & 2 & Yes & 0.9569 & Yes \\
ESDIRK32b & 4 & 0.2929 & 2  & Yes & 0 & Yes  & & 3 & Yes & 1.609 & Yes      \\
ESDIRK43a & 5 & 0.5728 & 4  & Yes & 0 & Yes  & & 3 & Yes & 0.5525 & Yes \\
ESDIRK43b & 5 & 0.4359 & 3  & Yes & 0 & Yes  & & 4 & Yes & 0.7175 & Yes  \\
ESDIRK54a & 7 & 0.26   & 5  & Yes & 0 & Yes  & & 4 & Yes & 0.7483 & Yes \\
ESDIRK54b & 7 & 0.27  & 4 & Yes & 0 & Yes  &   & 5 & Yes & 0.8732 & Yes  \\
\hline ESDIRK32c & 4 & 0.5  & 3  & Yes & 0 & Yes  &  & 2 & Yes & 1 & Yes \\
ESDIRK45c & 6 & 0.25 &  4 & Yes & 0 & Yes  & & (5) & No    & $\infty$ & No  \\
\hline
\end{tabular}
\end{center}
\end{table}

A family of ESDIRK methods suitable for integration of stiff
systems of differential equations as well as index-1 systems of
differential algebraic equations have been constructed. The
integration methods of order $p$ are A- and L-stable as well as
stiffly accurate. The embedded methods for error estimation are of
order $p+1$. They are neither A- nor L-stable. This is of little
concern, since local extrapolation is not applied, i.e. the next
step is computed using the basic integration method of order $p$.
The methods have $s = p+1$ stages, but the first stage is the same
as the last stage in the previous step (FSAL). Hence, the
effective number of stages in the methods is $s-1$. Methods have
been constructed for $p = \set{1,2,3}$. These methods are called
ESDIRK12, ESDIRK23 and ESDIRK34, respectively. In addition, the
methods are equipped with a continuous extension that satisfies
the order conditions of the basic integration method. Therefore,
the continuous extensions have the same order as the basic
integration methods.
%Finally, strategies for the step size
%selection, adaptive Jacobian evaluation, and iteration matrix
%factorization for the methods have been discussed.

%\section*{Acknowledgments}
%The author thanks 2-control ApS for financial support.

\bibliographystyle{siam}
\bibliography{jbjbib,qpbook}

\end{document}